\documentclass[12pt]{article}
\pdfoutput=1
\usepackage[utf8]{inputenc}
\usepackage{amsfonts,amsmath,amsthm,amssymb}
\usepackage{nicefrac}
\usepackage{enumitem}
\usepackage[colorlinks,allcolors=blue]{hyperref}
\usepackage{url}
\usepackage[longnamesfirst, authoryear]{natbib}
\usepackage[title]{appendix}

\newcommand{\red}{\color{black}}
\newcommand{\black}{\color{black}}

\newtheorem{theorem}{Theorem}

\newtheorem{lemma}{Lemma}

\newtheorem{proposition}{Proposition}

\newtheorem{corollary}{Corollary}

\theoremstyle{remark}
\newtheorem{remark}{Remark}

\newtheorem{example}{Example}

\newtheorem{assumption}{A}

\let\Pr\relax
\DeclareMathOperator\Pr{\mathbb{P}}
\DeclareMathOperator\Qr{\mathbb{Q}}
\DeclareMathOperator\E{\mathbb{E}}

\DeclareMathOperator\var{var}

\newcommand{\KL}[2]{ { D \left({#1} \;\middle\Vert\; {#2}\right) } }

\newcommand{\calF}{{\mathcal F}}

\newcommand{\calM}{{\mathcal M}}

\newcommand{\Intp}[1]{ { {\mathbb Z_{\ge 0}}^{#1} } }

\newcommand{\Real}[1]{ { {\mathbb R}^{#1} } }
\newcommand{\Realp}[1]{ { {\mathbb R}_{\ge 0}^{#1} } }

\newcommand{\dss}{\displaystyle}
\newcommand{\one}{{\mathbf  1}}
\newcommand{\inv}{^{-1}}

\newcommand{\tran}{^{\top}}

\title{\textbf{Optimal Chernoff and Hoeffding Bounds for Finite 
\red State \black 
Markov Chains}}
\author{
Vrettos Moulos
\thanks{Supported in part by the NSF grant CCF-1816861.}
\\
UC Berkeley \\
\href{mailto:vrettos@berkeley.edu}{vrettos@berkeley.edu}
\and
Venkat Anantharam
\thanks{Supported in part by the NSF grants 
\red CNS-1527846, CCF-1618145, CCF-1901004,
\black the NSF Science \& Technology Center grant CCF-0939370 (Science of Information), and the William and Flora Hewlett Foundation supported Center for Long Term Cybersecurity at Berkeley.}
\\
UC Berkeley \\
\href{mailto:ananth@berkeley.edu}{ananth@berkeley.edu}
}
\date{}

\begin{document}

\maketitle

\begin{abstract}
This paper develops an optimal Chernoff type bound for the probabilities
of large deviations of sums $\sum_{k=1}^n f (X_k)$ where
$f$ is a real-valued function and $(X_k)_{k \in \mathbb{Z}_{\ge 0}}$ is a finite 
state Markov chain with an arbitrary initial distribution and an irreducible 
transition probability
matrix 
satisfying a mild assumption on its positivity pattern, related to the function $f$ being considered.
The novelty lies in this being a non-asymptotic finite sample bound. Further,
our bound 
is optimal in the large deviations 
 sense, attaining a
constant prefactor and an exponential decay with the optimal large deviations rate.
Moreover, through
a Pinsker type inequality and a Hoeffding type lemma, we are able to loosen up our Chernoff type bound to a Hoeffding type bound and reveal the sub-Gaussian nature of the sums.
Finally,
under the same mild assumption on the positivity pattern
of the transition probability matrix,
we prove a uniform multiplicative ergodic theorem for 
the exponential family of tilted transition probability matrices corresponding to $f$.
\end{abstract}

\section{Introduction}

Let $S$ be a finite 
%state space 
\red set \black
and $(X_k)_{k \in \Intp{}}$ 
%be the 
\red the \black
coordinate process on $S^\Intp{}$.
Given an initial distribution $q$ on $S$, and a stochastic matrix
$P$, there exists a unique probability measure $\Pr_q$ on the sequence space such that the coordinate process $(X_k)_{k \in \Intp{}}$ is a Markov chain 
\red
with transition probability matrix $P$,
with respect to the filtration of $\sigma$-fields
$(\calF_n := \sigma (X_0, \ldots, X_n), n \ge 0)$.
\black
If we assume further that $P$ is irreducible, then there exists a unique stationary distribution $\pi$
\red for the transition probability matrix $P$, 
\black
and for any real-valued function $f : S \to \Real{}$ the empirical mean
$n\inv \sum_{k=1}^n f (X_k)$ converges $\Pr_q$-almost-surely
to the stationary mean $\pi (f) := \sum_x f (x) \pi (x)$.
The goal of this work is to quantify the rate of this convergence
by developing finite sample upper bounds for the large deviations probability
\[
\Pr_q \left(\frac{1}{n} \sum_{k=1}^n f (X_k) \ge \mu\right),
~\text{for}~ \mu \ge \pi (f).
\]

The significance of studying finite sample bounds for such tail probabilities is not only theoretical but also practical, since concentration inequalities for Markov dependent random variables have wide applicability in statistics, computer science and learning theory.
Just to mention a few applications, first and foremost this convergence forms the backbone behind all Markov chain Monte Carlo (MCMC) integration techniques, see~\cite{MCMC53}.
Moreover, tail bounds of this form have been used by~\cite{JSV01} to develop an approximation algorithm for the permanent of a nonnegative matrix.
In addition, in the stochastic multi-armed bandit literature the analysis of learning algorithms is based on tail bounds of this type,
see the survey of~\cite{BC12}.
More specifically the work of~\cite{Moulos19-bandits-identification} uses such a bound to tackle a Markovian identification problem.

\subsection{Chernoff Bound}

The classic large deviations theory for Markov chains due to~\cite{Miller61, Donsker-Varadhan-I-75, Gartner77, Ellis84, Dembo-Zeitouni-98} suggests that asymptotically the large deviations probability decays exponentially and the rate is given by the convex conjugate $\Lambda^* (\mu)$
of the log-Perron-Frobenius eigenvalue $\Lambda (\theta)$ of the nonnegative irreducible matrix $\tilde{P}_\theta (x, y) := P (x, y) e^{\theta f (y)}$. In particular
\[
\lim_{n \to \infty}
\frac{1}{n} \log
\Pr_q \left(\frac{1}{n} \sum_{k=1}^n f (X_k) \ge \mu\right) =
- \Lambda^* (\mu), ~ \text{for} ~ \mu \ge \pi (f).
\]
Our objective is to develop a finite sample bound which captures this exponential decay and has a constant prefactor that does not depend on $\mu$,
and is thus useful in applications. 
A counting based approach by~\cite{Davisson-Longo-Sgarro-81} is able to capture this exponential decay but with a suboptimal prefactor that depends polynomially on $n$.
Through the development in the book of~\cite{Dembo-Zeitouni-98} (Theorem 3.1.2), 
which is also presented by~\cite{WH17},
one is able to obtain a constant prefactor, which though depends on $\mu$.
This is unsatisfactory because exact large deviations for Markov chains, see~\cite{Miller61,KM03}, yield that,
at least when the supremum $\sup_{\theta \in \Real{}} \{\theta \mu - \Lambda (\theta)\} = \Lambda^* (\mu)$ is attained at $\theta_\mu$, then
\[
\Pr_q \left(\frac{1}{n} \sum_{k=1}^n f (X_k) \ge \mu\right) \sim 
\frac{\E_{X \sim q} [v_{\theta_\mu} (X)]}
{\theta_\mu \sqrt{2 \pi n \sigma_{\theta_\mu}^2}} e^{-n \Lambda^* (\mu)}, 
~\text{as}~ n \to \infty,
\]
where $\sigma_{\theta_\mu}^2 = \Lambda'' (\theta_\mu)$ and $v_{\theta_\mu}$
is a right Perron-Frobenius eigenvector of $\tilde{P}_{\theta_\mu}$.
Here $\sim$ denotes that the ratio of the expressions on the left hand side and the right hand side converges to $1$, and $\Lambda'' (\theta_\mu)$ denotes the second derivative in $\theta$ of $\Lambda (\theta)$ at $\theta = \theta_\mu$. 
Thus, if we allow dependence on $\mu$, then the prefactor should be able to capture a decay of the order $1/\sqrt{n}$.
If we insist on 
%a constant prefactor 
\red no dependence on $\mu$ 
\black
though, the best that we can hope for is a constant prefactor, because otherwise we will contradict the central limit theorem for Markov chains.
This is argued formally at the end of~\autoref{sec:Chernoff}.

In our work we establish a tail bound with the optimal rate of exponential decay and a constant prefactor which depends only on the function $f$ and the stochastic matrix $P$,
under the following conditions on $P$.
Let $a := \min_x\: f (x)$, and $b := \max_x\: f (x)$. \red Based on $f$, \black we define two set of states, 
%the ones of maximum $f$ value 
\red
$S_b := \{x \in S : f (x) = b\}$
\black
and 
%the ones of minimum $f$ value 
\red
$S_a := \{x \in S: f (x) = a\}$.
\black
We will require that $P$ satisfies some subset of the following structural assumptions
on the positivity pattern of $P$. 
%In particular we will
\red We will \black 
enforce~\autoref{as:1}-\autoref{as:2} for upper tail bounds,~\autoref{as:3}-\autoref{as:4} for lower tail bounds, and~\autoref{as:1}-\autoref{as:4} when we want to bound both tails.
\begin{assumption}\label{as:1}
The submatrix of $P$ with rows and columns in $S_b$ is irreducible.
\end{assumption}
\begin{assumption}\label{as:2}
For every $x \in S - S_b$, there exists 
$y \in S_b$ such that $P (x, y) > 0$.
\end{assumption}
\begin{assumption}\label{as:3}
The submatrix of $P$ with rows and columns in $S_a$ is irreducible.
\end{assumption}
\begin{assumption}\label{as:4}
For every $x \in S - S_a$, there exists 
$y \in S_a$ such that $P (x, y) > 0$.
\end{assumption}
\red As we will see shortly, with 
these \black
assumptions we are essentially enforcing that after suitable tilts of the transition probability matrix we are able to produce new Markov chains that can realize any stationary mean in $(a, b)$.
\red
Our assumptions are general enough to capture all Markov chains, reversible or not, for which all the transitions have a positive probability.
\black
%as well as all finitely supported IID sequences. 

The key technique to derive our Chernoff type bound is the old idea due to~\cite{Esscher32} of an exponential tilt, which lies at the heart of large deviations
\red theory. \black
In the world of statistics those exponential changes of measure go by the name exponential families and the standard reference is the book of~\cite{Brown86}.
Exponential tilts of stochastic matrices generalize those of finitely supported probability distributions, and were first introduced in the work of~\cite{Miller61}.
Subsequently they formed one of the main tools in the study of large deviations for Markov chains, see~\cite{Donsker-Varadhan-I-75,Gartner77,Ellis84,Dembo-Zeitouni-98,BM00,KM03}.
Naturally they are also the key object when one conditions on the 
%second-order empirical distribution 
\red pair empirical distribution \black
of a Markov chain and considers conditional limit 
\red theorems, as \black
in~\cite{Csizar-Cover-Choi-87,Bolthausen-Schmock-89}.
A more recent development by~\cite{Nagaoka-05}
gives an information geometry perspective to this concept,
while~\cite{HW16} examine the problem of parameter estimation for exponential families of stochastic matrices.

Here we build on exponential families of stochastic matrices and 
\red
by studying the analyticity properties of 
the Perron-Frobenius eigenvalue and its 
associated eigenvector as we parametrically move
the mean of $f$ under exponential tilts, 
together with conjugate duality,
\black
%together with some Perron-Frobenius theory,
%analyticity of Perron-Frobenius eigenvalues and eigenvectors, as well as conjugate duality 
we are able to establish our main Chernoff type bound.

\begin{theorem}\label{thm:Chernoff-bound}
Let $P$ be an irreducible stochastic matrix on the finite state space $S$, with stationary distribution $\pi$, 
\red which, combined \black with a real-valued function 
\red $f : S \to \Real{}$, satisfies~\autoref{as:1}-\autoref{as:2}. \black
\red Then, for any initial distribution $q$, we have 
\black
\[
\Pr_q \left(\frac{1}{n} \sum_{k=1}^n f (X_k) \ge \mu\right) \le K_u e^{-n \Lambda^* (\mu)},
~\text{for}~ \mu \ge \pi (f),
\]
where $K_u = K_u (P, f)$ is the constant from~\autoref{prop:evec-ratio}, and depends only on the stochastic matrix $P$ and the function $f$.
\end{theorem}
\begin{remark}
Since $f$ is arbitrary and our assumptions~\autoref{as:1}-\autoref{as:2} 
and~\autoref{as:3}-\autoref{as:4} are symmetric, we can substitute $f$ with $-f$, so that~\autoref{thm:Chernoff-bound} yields a Chernoff type bound for the lower tail as well. In particular, assuming~\autoref{as:3}-\autoref{as:4} we have
\[
\Pr_q \left(\frac{1}{n} \sum_{k=1}^n f (X_k) \le \mu\right) \le K_l e^{-n \Lambda^* (\mu)},
~\text{for}~ \mu \le \pi (f),
\]
where $K_l = K_u (P, -f)$.
\end{remark}
\begin{remark}
Similarly assuming~\autoref{as:1}-\autoref{as:4} we have the following two-sided Chernoff type bound.
\[
\Pr_q \left(\frac{1}{n} \sum_{k=1}^n f (X_k) \in F\right) \le
2 K e^{-n \inf_{\mu \in F} \Lambda^* (\mu)},
~\text{for any}~ F ~\text{closed in}~ \Real{},
\]
where $K = \max \{K_l, K_u\}$.
\end{remark}
\begin{remark}
According to~\autoref{prop:evec-ratio}, when $P$ is a positive stochastic matrix, i.e. all the transitions have positive probability, we can replace $K$ with
\[
K \le \dss \max_{x,y,z} \frac{P (x, z)}{P (y, z)}.
\]
\end{remark}
\begin{remark}
According to~\autoref{prop:evec-ratio}, when $P$ induces an IID sequence, i.e. all the rows of $P$ are identical, then $K = 1$.
Thus~\autoref{thm:Chernoff-bound} generalizes the classic bound of~\cite{Chernoff52} for finitely supported IID sequences.
\end{remark}

\subsection{Hoeffding Bound}

Although Chernoff type bounds for Markov chains have not been extensively studied in the literature, and that's exactly the focus of this work, there is a vast literature on Hoeffding type inequalities for Markov chains.
\cite{Gillman93} obtained the first Hoeffding type bound for \red reversible finite state Markov chains. \black
Reversibility is a key assumption in his work because it 
%leads 
\red allows one to restrict attention \black
to self-adjoint operators and then it is possible to apply the matrix perturbation theory of~\cite{Kato66} in order to derive a bound on the largest eigenvalue of the 
\red
%perturbed 
self-adjoint operator $\tilde{P}_\theta$ defined 
in~\eqref{eq:twist} below.
\black
Later on~\cite{Dinwoodie95} obtained an improved prefactor.
Using the same spectral techniques~\cite{Lezaud98} obtained a Bernstein type inequality which is also applicable to some nonreversible finite
\red state \black Markov chains, and
which was later improved in the work of~\cite{Paulin15}.
\cite{Kahale97} introduced the idea of reducing the problem 
\red to a two state chain, which turned out to be 
\black 
very fruitful. \cite{LP04} employed this idea 
\red and, by performing exact calculations, they obtained 
\black 
a bound which is optimal for two state chains in the large deviations sense, as well as a Hoeffding type bound with variance proxy
$\frac{1 + \lambda \lor 0}{1 - \lambda \lor 0} \frac{(b-a)^2}{4}$, where $\lambda$ is the second largest eigenvalue of the reversible stochastic matrix $P$, as opposed to the classic variance proxy
for IID sequences $\frac{(b-a)^2}{4}$ due to~\cite{Hoeffding63}. \cite{Miaso14} extended this work to general state spaces without the reversibility assumption,~\cite{Rao19}
considered 
%finite stationary 
\red stationary finite state \black
Markov chains but allowed time-varying functions $f_i$, %and finally~\cite{FJS18-Bernstein, FJS18-Hoeffding}
\red and finally~\cite{FJS18-Bernstein} and \cite{FJS18-Hoeffding} \black
obtained both Bernstein and Hoeffding type bounds for general state space Markov chains and time-varying functions $f_i$.

Here we develop a Hoeffding type bound by loosening up our Chernoff type bound in~\autoref{thm:Chernoff-bound} using a Pinsker type
inequality in~\autoref{lem:Pinsker}. In the process a Hoeffding type lemma, in~\autoref{lem:Hoeffding}, is established as the dual of our Pinsker type inequality.

\begin{theorem}\label{thm:Hoeffding-bound}
Let $P$ be an irreducible stochastic matrix on the finite state space $S$, with stationary distribution $\pi$, 
\red which, combined with \black a real-valued function 
\red $f : S \to [a, b]$, satisfies~\autoref{as:1}-\autoref{as:2}. \black
\red Then, for any initial distribution $q$, we have
\black
\[
\Pr_q \left(\frac{1}{n} \sum_{k=1}^n f (X_k) \ge \mu\right) \le
K_u e^{-n \frac{(\mu - \pi (f))^2}{2 \sigma_u^2}} \le 
K_u e^{-n \frac{2 (\mu - \pi (f))^2}{(b-a + 2 K_u L_u)^2}},
~\text{for}~ \mu \ge \pi (f),
\]
where $\sigma_u^2 = \sigma_u^2 (P, f) := \sup_{\theta \in \Realp{}}\: \Lambda'' (\theta) < \infty$, $\Lambda'' (\theta)$ denotes
the second derivative of $\Lambda (\theta)$ in $\theta$, and
$K_u = K_u (P, f), ~ L_u = L_u (P, f)$ are the constants from~\autoref{prop:evec-ratio}.
\end{theorem}
\begin{remark}
Since $f$ is arbitrary and our assumptions~\autoref{as:1}-\autoref{as:2} 
and~\autoref{as:3}-\autoref{as:4} are symmetric, we can substitute $f$ with $-f$, so that~\autoref{thm:Hoeffding-bound} yields a Hoeffding type bound for the lower tail as well. In particular, assuming~\autoref{as:3}-\autoref{as:4} we have
\[
\Pr_q \left(\frac{1}{n} \sum_{k=1}^n f (X_k) \le \mu\right) \le
K_l e^{-n \frac{(\mu - \pi (f))^2}{2 \sigma_l^2}} \le 
K_l e^{-n \frac{2 (\mu - \pi (f))^2}{(b-a + 2 K_l L_l)^2}},
~\text{for}~ \mu \le \pi (f),
\]
\red
where $\sigma_l^2 := \sigma_u^2 (P, -f)$,
and $L_l := L_u (P, -f)$.
\black
\end{remark}
\begin{remark}
According to~\autoref{prop:evec-ratio}, when $P$ induces an IID sequence, i.e. all the rows of $P$ are identical, then $K = \max \{K_u, K_l\} = 1$ and $L = \max \{L_u, L_l\} = 0$.
Thus~\autoref{thm:Hoeffding-bound} generalizes the classic bound of~\cite{Hoeffding63} for finitely supported IID sequences.
\end{remark}
\begin{remark}
Our variance proxy
$\sigma_u^2 = \sup_{\theta \in \Realp{}} \Lambda'' (\theta) \le (b-a + 2 K_u L_u)^2/4$,
according to~\autoref{lem:dual-map}, has an interpretation as a worst
case variance among all the tilted Markov chains,
and thus parallels the variance proxy from the IID case which is the supremum of the variances among the tilted \red distributions, and which can be \black upper bounded by $(b-a)^2/4$.
\end{remark}

\subsection{Organization of Paper}

The rest of the paper proceeds as follows.
\autoref{sec:exp-fam} contains the classic construction of exponential families of stochastic matrices, the duality between the canonical and mean parametrization, as well as many other useful properties for our bounds.
In~\autoref{sec:Chernoff} and~\autoref{sec:Hoeffding} we analyze the limiting behavior of the family under our assumptions~\autoref{as:1}-\autoref{as:2}, and we establish our Chernoff~(\autoref{thm:Chernoff-bound}) and Hoeffding~(\autoref{thm:Hoeffding-bound}) type bounds.
Finally in~\autoref{sec:ergodic} we develop a uniform multiplicative ergodic theorem~(\autoref{thm:ergodic}).

\section{Exponential Family of Stochastic Matrices}\label{sec:exp-fam}

\subsection{Construction}

Exponential tilting of stochastic matrices originates in the work of~\cite{Miller61}. 
\red Following this, we define \black an exponential family of stochastic matrices which is able to produce Markov chains with shifted stationary means.
The generator of the exponential family is an irreducible stochastic matrix $P$, 
which for this section is not assumed to satisfy~\autoref{as:1}-\autoref{as:4},
and $\theta \in \Real{}$ represents the canonical parameter of the family.
Then we define
\red
\begin{equation}        \label{eq:twist}
\tilde{P}_\theta (x, y):= P (x, y) e^{\theta f (y)},
\end{equation}
\black
(or $\widetilde{\left(P\right)}_\theta (x, y)$, where $\widetilde{(\cdot)}_\theta$ is thought as an operator over matrices). 
$\tilde{P}_\theta$ has the same nonnegativity structure as $P$,
hence it is irreducible and we can use the Perron-Frobenius theory in order to normalize it and turn it into a stochastic matrix.
Let $\rho (\theta)$ (or $\rho (\tilde{P}_\theta)$) be the spectral radius of $\tilde{P}_\theta$,
which from the Perron-Frobenius theory 
%is an simple eigenvalue 
\red is a simple eigenvalue \black
of $\tilde{P}_\theta$, called the Perron-Frobenius eigenvalue, associated with unique left and right eigenvectors
$u_\theta, ~ v_\theta$
(or $u_{\tilde{P}_\theta},~v_{\tilde{P}_\theta}$)
such that they 
%are both positive, 
\red both have all entries strictly positive, \black
$\sum_x u_\theta (x) = 1,$  and $\sum_x u_\theta (x) v_\theta (x) = 1$, see for instance Theorem 8.4.4 in the book of~\cite{HJ13}.
Using $\tilde{P}_\theta$ we define a family of nonnegative irreducible matrices, parametrized by $\theta$, in the following way
\red
\begin{equation}    \label{eq:expfamily}
\left(P\right)_\theta (x, y) = P_\theta (x, y) := \frac{\tilde{P}_\theta (x, y) v_\theta (y)}
{\rho (\theta) v_\theta (x)},
\end{equation}
\black
which are stochastic, since
\[
\sum_y P_\theta (x, y) = \frac{1}{\rho (\theta) v_\theta (x)} \cdot \sum_y \tilde{P}_\theta (x, y) v_\theta (y) = 1, ~ \text{for} ~ x \in S.
\]
In addition 
%their stationary distributions 
\red the stationary distributions of the $P_\theta$ \black
are given by
\[
\pi_\theta (x) := u_\theta (x) v_\theta (x), ~ \text{for} ~ x \in S,
\]
since
\red
\[
\sum_x \pi_\theta (x) P_\theta (x, y) =
\frac{v_\theta (y)}{\rho (\theta)} \cdot 
\sum_x u_\theta (x) \tilde{P}_\theta (x, y) =
\pi_\theta (y),
~ \text{for} ~ y \in S.
\]
\black
Note that the generator stochastic matrix, $P$, is the member of the family that corresponds to $\theta = 0$, i.e. $P_0 = P, ~ \rho (0) = 1, ~ u_0 = \pi, ~ v_0 = \one$, and $\pi_0 = \pi$, where $\one$ is the all ones vector. 
%In addition 
\red In general \black
it is possible that the family is degenerate as the following example suggests.
\begin{example}\label{ex:deg}
Let $S = \{\pm 1\}, ~ P (x, y) = 1 \{x \neq y\}$, and $f (x) = x$.
Then $\rho (\theta) = 1, ~ v_\theta (-1) = \frac{1+e^\theta}{2}, ~v_\theta (1) = \frac{1+e^{-\theta}}{2}$, and $P_\theta = P$ for any $\theta \in \Real{}$.
\end{example}
A basic property 
%of the family 
\red of the exponential family $P_\theta$ \black
is that the composition of $(\cdot)_{\theta_1}$ with $(\cdot)_{\theta_2}$, is the transform $(\cdot)_{\theta_1 + \theta_2}$, and so composition is commutative. Furthermore we can undo the transform $(\cdot)_{\theta}$ by applying $(\cdot)_{-\theta}$.
\red We state this formally for convenience. \black
\begin{lemma}\label{lem:comp}
For any irreducible stochastic matrix $P$, and any $\theta_1, \theta_2 \in \Real{}$
\[
\left(\left(P\right)_{\theta_2}\right)_{\theta_1} = 
\left(P\right)_{\theta_1 + \theta_2}.
\]
\begin{proof}
It suffices to check that
$\left(\frac{v_{\theta_1 + \theta2} (y)}{v_{\theta_2} (y)}, ~ y \in S\right)$
is a right eigenvector of the matrix with entries
$\left(\frac{P (x, y) e^{\theta_2 f (y)} v_{\theta_2} (y)}{\rho (\theta) v_{\theta_2} (x)}\right)e^{\theta_1 f (y)}$,
with the corresponding eigenvalue being $\frac{\rho (\theta_1 + \theta_2)}{\rho (\theta_2)}$. 
This is a straightforward calculation.
\end{proof}
\end{lemma}

\subsection{Mean Parametrization}

The exponential family 
\red $P_\theta$ defined in~\eqref{eq:expfamily} \black can be reparametrized using the mean parameters $\mu = \pi_\theta (f)$. The duality between the canonical parameters $\theta$ and the mean parameters $\mu$
is manifested through the log-Perron-Frobenius eigenvalue
$\Lambda (\theta) := \log \rho (\theta)$.
More specifically, from~\autoref{lem:dual-map} it follows that there are two cases for the mapping
$\theta \mapsto P_\theta$.
In the nondegenerate case that this mapping is nonconstant, $\Lambda' (\theta)$ is a strictly increasing bijection between the set $\Real{}$ of canonical parameters and the set
\red
\begin{equation}    \label{eq:parameterinterior}
\calM := \{\mu \in \Real{} : \pi_\theta (f) = \mu, ~ \text{for some} ~ \theta \in \Real{}\}
\end{equation}
\black 
of mean parameters, which is an open interval.
Therefore, with some abuse of notation, for any $\mu \in \calM$ we may write $u_\mu, v_\mu, P_\mu, \pi_\mu$ for $u_{{\Lambda'}\inv (\mu)}, v_{{\Lambda'}\inv (\mu)}, P_{{\Lambda'}\inv (\mu)}, \pi_{{\Lambda'}\inv (\mu)}$. In the degenerate case that the mapping is constant, $\Lambda' (\theta) = \pi (f)$, and the set $\calM$ is the singleton $\{\pi (f)\}$. An illustration of the degenerate case is~\autoref{ex:deg}.

\begin{lemma}\label{lem:dual-map}
Let $P$ be an irreducible stochastic matrix, and $f : S \to \Real{}$ a real-valued function on the state space $S$. Then
\begin{enumerate}[label=(\alph*)]
    \item $\rho (\theta), ~\Lambda (\theta), ~ u_\theta$ and $v_\theta$ are analytic functions of $\theta$ on $\Real{}$.
    \item $\Lambda' (\theta) = \pi_\theta (f)$.
    \item 
    $
    \Lambda'' (\theta) =
    \var_{(X, Y) \sim \pi_\theta \odot P_\theta}
    \left(
    f (Y) + \frac{v_\theta (X)}{v_\theta (Y)}
    \frac{d}{d \theta} \frac{v_\theta (Y)}{v_\theta (X)}
    \right)
    $, where $\pi_\theta \odot P_\theta$ denotes the bivariate distribution defined by
    $(\pi_\theta \odot P_\theta) (x, y) := \pi_\theta (x) P_\theta (x, y)$.
    \item Either $P_\theta = P_0 = P$ for all $\theta \in \Real{}$ (degenerate case), or $\theta \mapsto P_\theta$ is an injection (nondegenerate case).
    
    Moreover, in the degenerate case $\Lambda (\theta) = \pi_0 (f) \theta$ is linear, while in the nondegenerate case $\Lambda (\theta)$ is strictly convex.
\end{enumerate}
\end{lemma}
The proof of~\autoref{lem:dual-map} can be found in~\autoref{app:misc-proofs}.

\subsection{Relative Entropy Rate and Conjugate Duality}
For two probability distributions $\Qr$ and $\Pr$ over the same measurable space we define the \emph{relative entropy} between $\Qr$ and $\Pr$ as
%\[
%\KL{\Qr}{\Pr} :=
%\begin{cases}
%\E_{\Qr} \left[\log \frac{d \Qr}{d \Pr} \right], &\text{if}~ \Qr \ll \Pr, \\
%\infty, &\text{otherwise}.
%\end{cases}
%\]
\red
\[
\KL{\Qr}{\Pr} :=
\begin{cases}
\E_{\Qr} \left[\log \frac{d \Qr}{d \Pr} \right], &\text{if
$Q$ is absolutely continuous with respect to $P$,}\\
\infty, &\text{otherwise}.
\end{cases}
\]
\black
Relative entropies of stochastic processes are most of the time trivial, and so we resort to the notion of relative entropy rate.
Let $Q, P$ be two stochastic matrices over the same state space $S$.
We further assume that $Q$ is irreducible with associated stationary distribution $\pi_Q$. 
For any initial distribution $q$ on $S$ we define the \emph{relative entropy rate} between the Markov chain $\Qr_q$ induced by $Q$ with initial distribution $q$, and the Markov chain $\Pr_q$ induced by $P$ with initial distribution $q$ as
\[
\KL{Q}{P} :=
\lim_{n \to \infty} \frac{1}{n}
\KL{\Qr_q \mid_{\calF_n}}{\Pr_q \mid_{\calF_n}},
\]
where $\Qr_q \mid_{\calF_n}$ and $\Pr_q \mid_{\calF_n}$
denote the finite dimensional distributions of the probability measures restricted to the sigma algebra $\calF_n$.
\red Note that 
%indeed 
the definition \black
is independent of the initial distribution $q$, since we can easily see using ergodic theory that
\[
\KL{Q}{P} =
\sum_{x, y} \pi_Q (x) Q (x, y) \log \frac{Q (x, y)}{P (x, y)} =
\KL{\pi_Q \odot Q}{\pi_Q \odot P},
\]
where $\pi_Q \odot Q$ denotes the bivariate distribution
\[
(\pi_Q \odot Q) (x, y) := \pi_Q (x) Q (x, y),
\]
and we use the standard notational conventions $\log 0 = \infty, ~ \log \frac{\alpha}{0} = \infty ~ \text{if} ~ \alpha > 0$,
and $0 \log 0 = 0 \log \frac{0}{0} = 0$.

For stochastic matrices which are elements 
\red of the exponential family 
$P_\theta$ defined in~\eqref{eq:expfamily} \black 
we simplify the relative entropy rate notation as follows.
For $\theta_1, \theta_2 \in \Real{}$ and $\mu_1 = \Lambda' (\theta_1),~ \mu_2 = \Lambda' (\theta_2)$ we write
\[
\KL{\theta_1}{\theta_2}, \KL{\mu_1}{\mu_2}
:= \KL{\pi_{\theta_1} \odot P_{\theta_1}}
{\pi_{\theta_1} \odot P_{\theta_2}}.
\]
For those relative entropy rates~\autoref{lem:KL-rep} suggests an alternative representation based on the parametrization. Its proof can be found in~\autoref{app:misc-proofs}.
\begin{lemma}\label{lem:KL-rep}
Let $\theta_1, \theta_2 \in \Real{}$ and $\mu_1 = \Lambda' (\theta_1),~ \mu_2 = \Lambda' (\theta_2)$. Then
\[
\KL{\theta_1}{\theta_2} =
\Lambda (\theta_2) - \Lambda (\theta_1) - \mu_1 (\theta_2 - \theta_1).
\]
\end{lemma}
We further define the convex conjugate of $\Lambda (\theta)$ as
$\Lambda^* (\mu) := \sup_{\theta \in \Real{}}\: \{\theta \mu - \Lambda (\theta)\}$. 
Moreover, since we saw in~\autoref{lem:dual-map} that $\Lambda (\theta)$ is convex and analytic, we have that the biconjugate of $\Lambda (\theta)$ is $\Lambda (\theta)$ itself, i.e.
$\Lambda (\theta) = \sup_{\mu \in \Real{}}\: \{\mu \theta - \Lambda^* (\mu)\}$.
The convex conjugate $\Lambda^* (\mu)$ represents the rate of exponential decay for large deviation events, and in the following~\autoref{lem:conv-conj}, which is established in~\autoref{app:misc-proofs}, we derive a closed form expression for it.
\begin{lemma}\label{lem:conv-conj}
\red
\[
\Lambda^* (\mu) =
\begin{cases}
\KL{\mu}{\pi (f)},
&\text{if}~ \mu \in \calM, \\
\dss\lim_{\hat{\mu} \to \mu} \KL{\hat{\mu}}{\pi (f)},
&\text{if}~ \mu \in \partial \calM, \\
\infty, &\text{otherwise},
\end{cases}
\]
where $\calM$ is defined in~\eqref{eq:parameterinterior}.
\black
\end{lemma}
An inspection of how the supremum was obtained in the previous~\autoref{lem:conv-conj} yields the following~\autoref{cor:conv-conj}.
\begin{corollary}\label{cor:conv-conj}
\[
\Lambda^* (\mu) =
\begin{cases}
\dss \sup_{\theta \ge 0}\: \{\theta \mu - \Lambda (\theta)\},
& \text{if}~ \mu \ge \pi (f), \\
\dss \sup_{\theta \le 0}\: \{\theta \mu - \Lambda (\theta)\},
& \text{if}~ \mu \le \pi (f).
\end{cases}
\]
\end{corollary}
\section{Optimal Chernoff Bound}
\label{sec:Chernoff}

\red
\subsection{The Class of Stochastic Matrices of Interest}
\black

In order to develop our upper tail bounds we assume that the irreducible stochastic matrix $P$ satisfies~\autoref{as:1}-\autoref{as:2},
\red for the given function $f : S \to \mathbb{R}$. \black 
Under those conditions we are able to show in~\autoref{prop:evec-ratio} that the ratio of the entries of the right Perron-Frobenius eigenvector 
$v_\theta (y)/ v_\theta (x)$ is uniformly bounded for $\theta \in \Realp{}$.
%Moreover, 
%those conditions capture a large class of Markov chains, for instance 
\red Note that the conditions~\autoref{as:1}-\autoref{as:2} are satisfied by 
Markov chains where every transitions has a positive probability. For these Markov chains, and in particular for Markov chains that induce IID processes, 
%For those two categories 
we provide explicit uniform bounds in~\autoref{prop:evec-ratio}.
\black

The following example suggests that we cannot meet the requirement that the ratios of the entries of the right Perron-Frobenius eigenvector is uniformly bounded if we drop assumption~\autoref{as:1}.
\begin{example}
Let $S = \{\pm 1\}, ~ P (x, y) = 1 \{x=-1\}/2 + 1 \{x=1, y=-1\}$, and $f (x) = x$.
Then $\rho (\theta) =  \frac{1 + \sqrt{1 + 8 e^{2\theta}}}{4} e^{-\theta}$, and
$v_\theta (-1) / v_\theta (1) = \rho (\theta) e^\theta \to \infty$ as $\theta \to \infty$.
\end{example}
Similarly a birth-death chain 
%shows the necessity 
\red illustrates the role \black
of assumption~\autoref{as:2}.
\begin{example}
Let $S = \{-1, 0, 1\}, ~ P (x, y) = 1 \{x+y \neq 0\}/2$ and $f (x) = -x$. Then 
$\rho (\theta) = \frac{1}{4} e^{\theta} \left(1 + e^{-\theta} + e^{-2 \theta}  + \sqrt{1 + 2e^{-\theta} + -5 e^{-2 \theta} + 2 e^{-3 \theta} + e^{-4 \theta}} \right)$, and
$v_\theta (0) / v_\theta (1) = 
2 \rho (\theta) - e^{-\theta} \to \infty$ as $\theta \to \infty$.
\end{example}

The natural interpretation of 
%those conditions 
\red conditions~\autoref{as:1}-\autoref{as:2} \black
is that they allow us to create new Markov chains with any stationary mean in the interval $[\pi (f), b)$, by selecting appropriate tilting levels $\theta \in \Realp{}$.

\red
\subsection{The Limiting Behavior of the Exponential Family}
\black

Define the matrix
\[
\overline{P}_\theta (x, y) := 
e^{- \theta b} \tilde{P}_\theta (x, y) =
e^{- \theta b} P (x, y) e^{\theta f (y)},
\]
and note that $\rho (\overline{P}_\theta) =
e^{- \theta b} \rho (\theta)$, as well as
$u_{\overline{P}_\theta} = u_\theta, ~ v_{\overline{P}_\theta} = v_\theta$. Hence $\overline{P}_\theta$ will help us study the asymptotic behavior of $P_\theta$, since
\[
P_\theta (x, y) = 
\frac{\overline{P}_\theta (x, y) v_\theta (y)}
{\rho (\overline{P}_\theta) v_\theta (x)}.
\]
Note that
\[
\overline{P}_\infty (x, y) := 
\lim_{\theta \to \infty} \overline{P}_\theta (x, y) =
\begin{cases}
P (x, y), & \text{if} ~ y \in S_b, \\
0, & \text{otherwise}.
\end{cases}
\]
Due to the structure imposed on $P$ through~\autoref{as:1}-\autoref{as:2},
the following~\autoref{lem:PF}, 
which constitutes a simple extension of the Perron-Frobenius theory for matrices which are not necessarily irreducible, suggests that
$\rho (\overline{P}_\infty) > 0$ is a simple eigenvalue of 
$\overline{P}_\infty$, which is associated with unique left and right eigenvectors $u_\infty, v_\infty$ such that $u_\infty (x) > 0$ for $x \in S_b$ and $u_\infty (x) = 0$ for $x \not\in S_b$, $v_\infty$ is positive, $\sum_x u_\infty (x) = 1$ and $\sum_x u_\infty (x) v_\infty (x) = 1$.
\begin{lemma}\label{lem:PF}
Let $M \in \Realp{s \times s}$ be a nonnegative matrix
such that after a consistent renumbering of its rows and columns 
\red
we can assume, 
for some $k \in \{1, \ldots, s\}$,
that $M$ consists of
an irreducible square block $A \in \Realp{k \times k}$,
and a rectangular block $B \in \Realp{(s-k) \times k}$
such that none of the rows of $B$ is zero,
\black
assembled together in the following way
\[
M =
\begin{bmatrix}
A & 0 \\
B & 0
\end{bmatrix}.
\]
Then, $\rho (M) = \rho (A) > 0$ is a simple eigenvalue of $M$,
which we call the Perron-Frobenius eigenvalue, and is associated with unique left and right eigenvectors $u_M, v_M$ such that $u_M$ has its first $k$ coordinates positive
and its last $s-k$ coordinates equal to zero,
$v_M$ is positive, $\sum_{x=1}^k u_M (x) = 1$, and
$\sum_{x=1}^k u_M (x) v_M (x) = 1$.
\begin{proof}
Let $u_A, v_A$ be the unique left and right eigenvectors of $A$ corresponding to the Perron-Frobenius eigenvalue $\rho (A)$,
such that both of them are positive, $\sum_{x=1}^k u_A (x) = 1$ and $\sum_{x=1}^k u_A (x) v_A (x) = 1$. Observe that the vectors
\[
u_M =
\begin{bmatrix}
u_A \\
0
\end{bmatrix}
,~\text{and}~
v_M =
\begin{bmatrix}
v_A \\
B v_A / \rho (A)
\end{bmatrix},
\]
are left and right eigenvectors of $M$ with associated eigenvalue $\rho (A)$, and satisfy all the conditions.

In addition, any eigentriple
$\lambda,
~\begin{bmatrix} u_h\tran & u_l\tran\end{bmatrix}\tran,
~\begin{bmatrix} v_h\tran & v_l\tran\end{bmatrix}\tran$
of eigenvalue and corresponding left and right eigenvectors of $M$,
will certainly have $u_l = 0$, and gives rise to an eigentriple 
$\lambda, u_h, v_h$ for $A$. 
Therefore, $\rho (M) = \rho (A)$ and the uniqueness of $u_M, ~ v_M$
follows from the uniqueness of $u_A, v_A$.
\end{proof}
\end{lemma}
Note that from~\autoref{lem:PF} for $k=s$ we recover the classic Perron-Frobenius theorem
\red(which we have of course used 
in the proof of~\autoref{lem:PF}).
\black

A continuity argument for simple eigenvalues and their corresponding eigenvectors, enables us to describe the asymptotic behavior of $P_\theta$ in~\autoref{lem:limit}.
\begin{lemma}\label{lem:limit}
$(u_\theta, ~\rho (\overline{P}_\theta), ~ v_\theta) \to 
(u_\infty, ~\rho (\overline{P}_\infty), ~ v_\infty)$, as $\theta \to \infty$, and so the following is a well defined stochastic matrix
\[
P_\infty (x, y) :=
\lim_{\theta \to \infty} P_\theta (x, y) =
\frac{\overline{P}_\infty (x, y) v_\infty (y)}{\rho (\overline{P}_\infty) v_\infty (x)}.
\]
\begin{proof}
Note that $\overline{P}_\infty$
possess the structure of~\autoref{lem:PF}.
Consider~\autoref{lem:anal} in~\autoref{app:anal},
with $M$ taken to be $\overline{P}_\infty$.
For $W$ in a sufficiently small neighborhood of $M$ the function $g (W)$ identified in the proof of that lemma is analytic and equals
$\begin{bmatrix} u_W\tran & \rho (W) & v_W\tran \end{bmatrix}\tran$
for all $W$ in that neighborhood that have the structure in~\autoref{lem:PF}.
Now, since $\overline{P}_\theta \to \overline{P}_\infty$
as $\theta \to \infty$, we have $\overline{P}_\theta$ 
is in this neighborhood for all sufficiently large $\theta$,
and $\overline{P}_\theta$, being irreducible, satisfies the conditions of~\autoref{lem:PF}.
The conclusion is now immediate.
\end{proof}
\end{lemma}
\begin{remark}\label{rem:any-mean}
The combination of the extended Perron-Frobenius %theory~\autoref{lem:PF}, 
\red theorem in~\autoref{lem:PF} \black
and the limiting behavior of the exponential family
\red establised in~\autoref{lem:limit} \black imply that
\[
\pi_\theta (f) \to b ~\text{as}~ \theta \to \infty,
\]
which together with~\autoref{lem:dual-map} (b) means that any mean $\mu$ in the interval $[\pi (f), b)$ can be realized by some exponential tilt $\theta \in \Realp{}$.
\end{remark}

A critical ingredient to obtain our tail bounds is the following~\autoref{prop:evec-ratio} which states that under the assumptions~\autoref{as:1}-\autoref{as:2} the ratio of the entries of the right Perron-Frobenius eigenvector stays uniformly bounded.
\begin{proposition}\label{prop:evec-ratio}
Let $P$ be an irreducible stochastic matrix on the finite state space $S$, 
\red which, combined with \black 
a real-valued function 
\red $f : S \to \Real{}$, satisfies~\autoref{as:1}-\autoref{as:2}. \black
Then
\[
K_u := \sup_{\theta \in \Realp{}, x, y \in S} \frac{v_{\overline{P}_\theta} (x)}{v_{\overline{P}_\theta} (y)} < \infty,\quad \text{and}~
L_u := \sup_{\theta \in \Realp{}, x, y \in S}
\left|\frac{d}{d \theta} \frac{v_{\overline{P}_\theta} (x)}{v_{\overline{P}_\theta} (y)}\right| < \infty,
\]
where $K_u = K_u (P, f)$ and $L_u = L_u (P, f)$ are constants depending on the stochastic matrix $P$,
and the function $f$. In particular
\begin{itemize}
    \item if $P$ induces an IID process, i.e. $P$ has identical rows, then $K_u = 1$ and $L_u = 0$;
    
    \item if $P$ is a positive stochastic matrix, then
    $
    K_u \le \dss \max_{x,y,z} \frac{P (x, z)}{P (y, z)}.
    $
\end{itemize}
\begin{proof}
\autoref{lem:dual-map} yields that $\theta \mapsto v_\theta (x)/v_\theta (y)$ is continuous, and so in conjunction with~\autoref{lem:limit} we have that the ratio of the entries of the right Perron-Frobenius eigenvector is uniformly bounded, hence $K_u < \infty$.

\red Moreover, using \black the chain rule we see 
\red that, for $\theta > 0$, we have \black
\[
\frac{d}{d \theta} \frac{v_{\overline{P}_\theta} (x)}{v_{\overline{P}_\theta} (y)} =
\sum_{z, w : P (z, w) > 0} 
- (b - f (w))
e^{-\theta (b - f (w))}P (z, w)
\frac{\partial}{\partial W (z, w)} \frac{v_W (x)}{v_W (y)} \mid_{W = \overline{P}_\theta}.
\]
To see why this formula holds, first observe that
$\overline{P}_\theta$,
being irreducible, satisfies the conditions of~\autoref{lem:PF}.
Next, observe that the last $s$ coordinates of $g (\overline{P}_\theta)$, in the notation of the proof
of~\autoref{lem:anal}, are all strictly positive.
With some abuse of notation since we are
not really thinking of $S$ as being enumerated, let us write the last $s$ coordinates of $g (W)$ as $g_W (s+1+x)$, for $x \in S$.~\autoref{lem:anal} then implies, that for all $x, y \in S$, the ratio $\frac{g_W (s+1+x)}{g_W (s+1+y)}$ 
is analytic in a sufficiently small neighborhood of $\overline{P}_\theta$.
Since a small enough variation in $\theta$ centered
around the given $\theta$ results in a variation of $W$
centered around $\overline{P}_\theta$ 
that lies in the set of matrices in this neighborhood that
satisfy~\autoref{lem:PF}
(in fact all such matrices are irreducible and, even further, are of the form $\overline{P}_{\theta'}$, for some $\theta'$),
we may use the notation $\frac{v_W (x)}{v_W (y)}$ for the ratio $\frac{g_W (s+1+x)}{g_W (s+1+y)}$.
The point is that what is really intended in the partial derivatives on the right hand side of the preceding equation is $\frac{\partial}{\partial W (z, w)} \frac{g_W (s+1+x)}{g_W (s+1+y)} \mid_{W = \overline{P}_\theta}$.

We claim that $\lim_{\theta \to \infty} \frac{d}{d \theta} \frac{v_{\overline{P}_\theta} (x)}{v_{\overline{P}_\theta} (y)} = 0$.
This is true because~\autoref{lem:anal} in~\autoref{app:anal} ensures that $\frac{\partial}{\partial W (z, w)} \frac{v_W (x)}{v_W (y)}$ is continuous at $\overline{P}_\infty$, more precisely
\[
\lim_{\theta \to \infty} \frac{\partial}{\partial W (z, w)} \frac{v_W (x)}{v_W (y)} \mid_{W = \overline{P}_\theta} =
\frac{\partial}{\partial W (z, w)} \frac{g_W (s+1+x)}{g_W (s+1+y)} \mid_{W = \overline{P}_\infty} \in \Real{}.
\]
Here, to be able to write the expression on the right hand
side of the preceding equation, we first observe that $\overline{P}_\infty$ satisfies the conditions
of~\autoref{lem:PF} and so the last $s$ coordinates of
$g (\overline{P}_\infty)$, in the notation of the proof
of~\autoref{lem:anal}, are all strictly positive, and so, by~\autoref{lem:anal}, for all $x, y \in S$ the ratio
$\frac{g_W (s+1+x)}{g_W (s+1+y)}$ is analytic in a neighborhood of $\overline{P}_\infty$. Further, the equality in the preceding equation is justified by the fact that, for all $\theta > 0$, 
$\frac{\partial}{\partial W (z, w)} \frac{v_W (x)}{v_W (y)} \mid_{W = \overline{P}_\theta}$
is just an alternate notation for
$\frac{\partial}{\partial W (z, w)} \frac{g_W (s+1+x)}{g_W (s+1+y)} \mid_{W = \overline{P}_\theta}$, and, 
for all $\theta$ large enough, $\overline{P}_\theta$
lies in the neighborhood around $\overline{P}_\infty$
guaranteed by~\autoref{lem:anal}.

Furthermore for the two cases for which we have a special handle on $K$ we argue as follows.
\begin{itemize}
\item Let $p$ be the probability distribution driving the IID process, i.e. all the rows of $P$ are identical and equal to $p$. Then we can see that 
$u_\theta (y) = \frac{p (y) e^{\theta f (y)}}{\sum_x p (x) e^{\theta f (x)}}, ~ \rho (\tilde{P}_\theta) = \sum_x p (x) e^{\theta f (x)},$ and $v_\theta (y) = 1$ for all $y \in S$,
since $\tilde{P}_\theta$ is the rank one matrix 
$\tilde{P}_\theta (x, y) = \rho (\tilde{P}_\theta) v_\theta (x) u_\theta (y)$.
    
\item If $P$ is a positive stochastic matrix then,
for any $\theta \in \Real{}, ~ x, y \in S$ we have that
\[
\frac{v_\theta (x)}{v_\theta (y)} =
\frac{\sum_z \tilde{P}_\theta (x, z) v_\theta (z)}{\sum_z \tilde{P}_\theta (y, z) v_\theta (z)} \le 
\max_{x,y,z} \frac{\tilde{P}_\theta (x, z)}{\tilde{P}_\theta (y, z)} = \max_{x,y,z} \frac{P (x, z)}{P (y, z)}.
\]
    
\end{itemize}
\end{proof}
\end{proposition}
Moreover under conditions~\autoref{as:1}-\autoref{as:2}
we are able to establish an explicit formula for the limiting
relative entropy rate.
\begin{lemma}\label{lem:limit-KL}
For any $\theta_2 \in \Real{}$, let $\mu_2 = \Lambda' (\theta_2)$. Then
\[ 
\KL{b}{\mu_2} :=
\lim_{\theta_1 \to \infty} \KL{\theta_1}{\theta_2} =
- \log \rho (\overline{P}_\infty) - (\theta_2 b - \Lambda (\theta_2)).
\]
\begin{proof}
From~\autoref{rem:any-mean} we have that
$\lim_{\theta \to \infty} \Lambda' (\theta) = b$, so
from~\autoref{lem:KL-rep} it suffices to show that
\[
\theta \Lambda' (\theta) - \Lambda (\theta) \to
- \log \rho (\overline{P}_\infty),
~\text{as}~ \theta \to \infty.
\]
Let $c = \max_{x \not \in S_b} f (x)$. 
Fix $x \in S$ and $y \not\in S_b$.
Pick $y_b \in S_b$ such that $P (x, y_b) > 0$ and as large as possible.
From~\autoref{prop:evec-ratio} we have that for any $\theta \in \Realp{}$
\[
P_\theta (x, y)
\le K_u \frac{P (x, y)}{P (x, y_b)} e^{- \theta (b-f (y))} P_\theta (x, y_b)
\le K_u K_0 e^{- \theta (b-c)},
\]
where
\[
K_0 := \max_{x \in S} \min_{y_b \in S} \frac{1}{P (x, y_b)} <
\infty.
\]
Therefore the stationary probability of any such $y$ is at most
$\pi_\theta (y) \le K_u K_0 e^{- \theta (b-c)}$, and so
\[
\pi_\theta (f) \ge (1- K_u K_0 |S| e^{- \theta (b-c)}) b + 
K K_0 |S| e^{- \theta (b-c)} a.
\]
From this we obtain that
$\lim_{\theta \to \infty} \theta (b - \Lambda' (\theta)) = 0$,
and the conclusion follows since from~\autoref{lem:limit} we have that
$\lim_{\theta \to \infty} (- \theta b + \Lambda (\theta)) = \log \rho (\overline{P}_\infty)$.

\end{proof}
\end{lemma}

\subsection{Chernoff Bound}

\begin{proof}[Proof of~\autoref{thm:Chernoff-bound}]
In order to derive our bounds we use a change of measure argument, an idea due to~\cite{Esscher32}.
We denote by $\Pr_q^{(\theta)}$
the probability distribution of the Markov chain with initial distribution $q$ and stochastic matrix $P_\theta$,
while for $\theta = 0$ we just write $\Pr_q$ for $\Pr_q^{(0)}$.
The finite dimensional distributions $\Pr_q \mid_{\calF_n}$ and $\Pr_q^{(\theta)} \mid_{\calF_n}$ are absolutely continuous with each other and their Radon-Nikodym derivative is given by
\[
\frac
{d \Pr_q \mid_{\calF_n}}
{d \Pr_q^{(\theta)} \mid_{\calF_n}} =
\frac
{v_{\theta} (X_0)}
{v_{\theta} (X_n)}
\exp\left\{
-\theta S_n + n  \Lambda (\theta)
\right\},
\]
where we denote the sums by $S_n := \sum_{k=1}^n f (X_k)$.

Fix $\theta \in \Realp{}$. Then
\begin{align*}
\Pr_q (S_n \ge n \mu) &=
\E_q \left[ 1 \{S_n \ge n \mu\} \right] \\
&= \E_q^{(\theta)} 
\left[
\frac
{v_{\theta} (X_0)}
{v_{\theta} (X_n)}
e^{-\theta S_n + n \Lambda (\theta)}
1 \{S_n \ge n \mu\}
\right] \\
&\le K_u 
\E_q^{(\theta)} 
\left[
e^{-\theta (S_n - n \mu)}
1 \{S_n \ge n \mu\}
\right]
e^{-n (\theta \mu - \Lambda (\theta))} \\
&\le K_u e^{-n (\theta \mu - \Lambda (\theta))},
\end{align*}
where in the first inequality we used~\autoref{prop:evec-ratio}.

When $\mu \in [\pi (f), b)$, we can set $\theta = {\Lambda'}^{-1} (\mu) \ge {\Lambda'}^{-1} (\pi (f)) = 0$
and then from~\autoref{lem:KL-rep} we have that
$\KL{\mu}{\pi (f)} = \theta \mu - \Lambda (\theta)$.
When $\mu = b$, we let $\theta$ go to $\infty$ and use~\autoref{lem:limit-KL}.
The conclusion follows from~\autoref{cor:conv-conj}.
\end{proof}

In this bound we cannot hope for something more than a constant prefactor. 
First of all, by differentiating twice the formula proved in~\autoref{lem:KL-rep} we obtain
\[
\lim_{\mu \to \pi (f)} \frac{1}{(\mu - \pi (f))^2} \KL{\mu}{\pi (f)} =
\frac{1}{2} \frac{1}{\Lambda'' (0)}.
\]

\red In addition, if we fix \black $z \ge 0$ and set $\mu = \pi (f) + c z/\sqrt{n}$, where $c^2 = \pi \left(\hat{f}^2 - (P \hat{f})^2\right)$ and $\hat{f}$ is a solution of the \emph{Poisson equation} $(I - P) \hat{f} = f - \pi (f)$,
then due to the central limit theorem for Markov chains, see for instance~\cite{Chung60}, we have that
\[
\lim_{n \to \infty} \Pr_q (S_n \ge n \mu) =
\int_z^\infty \frac{1}{\sqrt{2 \pi}} e^{-x^2/2} dx.
\]
Therefore if we want the optimal rate of exponential decay
and a prefactor which does not depend on $\mu$, then the best we can attain is a constant prefactor.
\section{Hoeffding Bound}\label{sec:Hoeffding}

In this section we relax the Chernoff type bound
to a Hoeffding type bound for finite Markov chains.
We achieve this by establishing a Pinsker type inequality in~\autoref{lem:Pinsker} for the exponential family of stochastic matrices which provides
a quadratic lower bound to the relative entropy rate.
Through conjugate duality this leads to a Hoeffding type lemma in~\autoref{lem:Hoeffding} which constitutes a quadratic
upper bound on the log-Perron-Frobenius eigenvalue.

We first develop a Pinsker type inequality for the exponential family of stochastic matrices under consideration.
\begin{lemma}\label{lem:Pinsker}
Let $P$ be an irreducible stochastic matrix on the finite state space $S$, which combined with a real-valued function $f : S \to [a, b]$ satisfies~\autoref{as:1}-\autoref{as:2}.
Then
\[
\KL{\mu}{\pi (f)} \ge
\frac{1}{2} \frac{(\mu - \pi (f))^2}{\sigma_u^2} \ge 
2 \frac{(\mu - \pi (f))^2}{(b-a + 2 K_u L_u)^2},
~\text{for any}~ \mu \ge \pi (f),
\]
where $\sigma_u^2 = \sigma_u^2 (P, f) := \sup_{\theta \in \Realp{}}\: \Lambda'' (\theta) \in (0, \infty)$,
and $K_u, L_u$ are the constants from~\autoref{prop:evec-ratio}.
\begin{proof}

We have that $\Lambda'' (\theta) > 0$ up to a null set, and so $\sigma^2 > 0$. \red Furthermore, from~\autoref{lem:dual-map} we have \black
\[
\Lambda'' (\theta) =
\var_{(X, Y) \sim \pi_\theta \odot P_\theta}
    \left(
    f (Y) + \frac{v_\theta (X)}{v_\theta (Y)}
    \frac{d}{d \theta} \frac{v_\theta (Y)}{v_\theta (X)}
    \right),
\]
\red
which, according to~\autoref{prop:evec-ratio}, is the \black variance of a random variable living in the interval $[a - K_u L_u, b + K_u L_u]$, and hence $\sigma^2$ is finite and upper bounded by $(b-a+2 K_u L_u)^2/4$.

Differentiating twice the formula proved in~\autoref{lem:KL-rep} we obtain
\[
\frac{\partial^2}{\partial \mu^2} \KL{\mu}{\pi (f)} =
\frac{1}{\Lambda'' ({\Lambda'}\inv (\mu))} \ge 
\frac{1}{\sigma_u^2},
~\text{for almost every}~ \mu \in \calM \cap [\pi (f), b).
\]
Using the fact that $\frac{\partial}{\partial \mu} \KL{\mu}{\pi (f)} \mid_{\mu = \pi (f)} = 0$, we conclude that
\[
\KL{\mu}{\pi (f)} = 
\int_{\pi (f)}^\mu (\mu - \nu)
\frac{\partial^2}{\partial \nu^2} \KL{\nu}{\pi (f)} d \nu \ge 
\frac{1}{2} \frac{(\mu - \pi (f))^2}{\sigma_u^2}.
\]

\end{proof}
\end{lemma}
Combining our \red Pinkser type inequality in~\autoref{lem:Pinsker}
and the Chernoff type bound in~\autoref{thm:Chernoff-bound},
\black
the Hoeffding type \red bound in~\autoref{thm:Hoeffding-bound} \black follows directly.

\begin{remark}\label{rem:Hoeffding}
Under \red assumptions~\autoref{as:3}-\autoref{as:4}, due to symmetry, we have \black that
\[
\KL{\mu}{\pi (f)} \ge
\frac{1}{2} \frac{(\mu - \pi (f))^2}{\sigma_l^2} \ge 
2 \frac{(\mu - \pi (f))^2}{(b-a + 2 K_l L_l)^2},
~\text{for any}~ \mu \le \pi (f),
\]
where $\sigma_l^2 = \sigma_u^2 (P, -f), ~K_l = K_u (P, -f)$,
and $L_l = L_u (P, -f)$.
\end{remark}

It is also possible to establish this Hoeffding bound directly
using the following Hoeffding type lemma for Markov chains
which is essentially the dual of our Pinsker type inequality.
\begin{lemma}\label{lem:Hoeffding}
Let $P$ be an irreducible stochastic matrix on the finite state space $S$, with stationary distribution $\pi$, which combined with a real-valued function $f : S \to [a, b]$ satisfies~\autoref{as:1}-\autoref{as:4}.
Then
\[
\Lambda (\theta) \le
\pi (f) \theta + \frac{1}{2} \sigma^2 \theta^2 \le 
\pi (f) \theta + \frac{(b-a + 2 K L)^2}{8} \theta^2,
~\text{for any}~ \theta \in \Real{},
\]
where $\sigma^2 = \max \{\sigma_u^2, \sigma_l^2\}, ~ K = \max \{K_u, K_l\}$,
and $L = \max \{L_u, L_l\}$.
\begin{proof}
%We plug 
\red Plugging in \black $\mu_2 = \pi (f)$ in~\autoref{lem:Pinsker},
and using~\autoref{lem:conv-conj} it is easy to see that
\[
\Lambda^* (\mu) \ge 
\frac{1}{2} \frac{(\mu - \pi (f))}{\sigma^2},
~\text{for all}~ \mu \in \Real{}.
\]
\red Finally, using \black the fact that $\Lambda (\theta)$ is the convex conjugate of $\Lambda^* (\mu)$ we conclude that
\[
\Lambda (\theta) \le \sup_{\mu \in \Real{}}\: 
\left\{\mu \theta - \frac{1}{2} \frac{(\mu - \pi (f))^2}{\sigma^2}\right\} =
\pi (f) \theta + \frac{1}{2} \sigma^2 \theta^2.
\]
\end{proof}
\end{lemma}
We call~\autoref{lem:Hoeffding} a Hoeffding type lemma,
since we will establish in~\autoref{sec:ergodic} that
$\frac{1}{n} \log
\E_q \left[
\exp\left\{\theta \sum_{k=1}^n f (X_k)\right\}
\right]$ converges to $\Lambda (\theta)$, and so at least asymptotically this reveals the sub-Gaussian structure of the sums.
\red
\section{A Uniform Multiplicative Ergodic Theorem}
\label{sec:ergodic}
\black

The classic linear ergodic theory for Markov chains,~\cite{Chung60} suggests that
\[
\frac{1}{n} \E_q \left[\sum_{k=1}^n f (X_k)\right]
\to \pi (f), ~\text{as}~ n \to \infty.
\]
\cite{BM00} and~\cite{KM03} have proved a multiplicative version of this under appropriate assumptions, which states that the \emph{scaled log-moment-generating-function} $\Lambda_n (\theta)$ converges pointwise to the log-Perron-Frobenius eigenvalue
\[
\Lambda_n (\theta) \to \Lambda (\theta), ~\text{as}~n \to \infty,
~\text{for any}~ \theta \in \Real{},
\]
where
\[
\Lambda_n (\theta) :=
\frac{1}{n} \log
\E_q \left[
\exp\left\{\theta \sum_{k=1}^n f (X_k)\right\}
\right].
\]
For our class of finite Markov chains we are able to establish a uniform multiplicative ergodic theorem in the terminology of~\cite{BM00}.
\begin{theorem}\label{thm:ergodic}
Let $P$ be an irreducible stochastic matrix on the finite state space $S$, which combined with a real-valued function $f : S \to \Real{}$ satisfies~\autoref{as:1}-\autoref{as:4}.
Then
\[
\sup_{\theta \in \Real{}}\: 
|\Lambda_n (\theta) - \Lambda (\theta)| \le 
\frac{\log K}{n}.
\]
where $K$ is the constant from~\autoref{prop:evec-ratio}.

Therefore $\Lambda_n (\theta)$ converges uniformly on $\Real{}$ to $\Lambda (\theta)$ as $n \to \infty$.
\begin{proof}
We start with the calculation
\begin{align*}
e^{n \Lambda_n (\theta)} 
&= \sum_{x_0, x_1, \ldots, x_{n-1}, x_n} q (x_0) P (x_0, x_1) e^{\theta f (x_1)} \cdots 
P (x_{n-1}, x_n) e^{\theta f (x_n)} \\
&= \sum_{x_0, x_n} q (x_0) \tilde{P}_\theta^n (x_0, x_n).
\end{align*}
From this, using the fact that $v_\theta$ is a right Perron-Frobenius eigenvector of $\tilde{P}_\theta$, we obtain
\[
\min_{x, y} \frac{v_\theta (y)}{v_\theta (x)}
\le \exp\left\{n \Lambda_n (\theta) - n \Lambda (\theta)\right\}
\le \max_{x, y} \frac{v_\theta (y)}{v_\theta (x)}.
\]
The conclusion now follows by applying~\autoref{prop:evec-ratio}.
\end{proof}
\end{theorem}

\section*{Acknowledgements} Vrettos Moulos would like to thank Jim Pitman and Satish Rao for many helpful discussions.

\bibliographystyle{apalike}
\bibliography{references}

\begin{thebibliography}{}

\bibitem[Balaji and Meyn, 2000]{BM00}
Balaji, S. and Meyn, S.~P. (2000).
\newblock Multiplicative ergodicity and large deviations for an irreducible
  {M}arkov chain.
\newblock {\em Stochastic Process. Appl.}, 90(1):123--144.

\bibitem[Bolthausen and Schmock, 1989]{Bolthausen-Schmock-89}
Bolthausen, E. and Schmock, U. (1989).
\newblock On the maximum entropy principle for uniformly ergodic {M}arkov
  chains.
\newblock {\em Stochastic Process. Appl.}, 33(1):1--27.

\bibitem[Brown, 1986]{Brown86}
Brown, L.~D. (1986).
\newblock {\em Fundamentals of statistical exponential families with
  applications in statistical decision theory}, volume~9 of {\em Institute of
  Mathematical Statistics Lecture Notes---Monograph Series}.
\newblock Institute of Mathematical Statistics, Hayward, CA.

\bibitem[Bubeck and Cesa-Bianchi, 2012]{BC12}
Bubeck, S. and Cesa-Bianchi, N. (2012).
\newblock {Regret Analysis of Stochastic and Nonstochastic Multi-armed Bandit
  Problems}.
\newblock {\em Foundations and Trends in Machine Learning}, 5(1):1--122.

\bibitem[Chernoff, 1952]{Chernoff52}
Chernoff, H. (1952).
\newblock A measure of asymptotic efficiency for tests of a hypothesis based on
  the sum of observations.
\newblock {\em Ann. Math. Statistics}, 23:493--507.

\bibitem[Chung, 1960]{Chung60}
Chung, K.~L. (1960).
\newblock {\em Markov chains with stationary transition probabilities}.
\newblock Die Grundlehren der mathematischen Wissenschaften, Bd. 104.
  Springer-Verlag, Berlin-G\"{o}ttingen-Heidelberg.

\bibitem[Csisz\'{a}r et~al., 1987]{Csizar-Cover-Choi-87}
Csisz\'{a}r, I., Cover, T.~M., and Choi, B.~S. (1987).
\newblock Conditional limit theorems under {M}arkov conditioning.
\newblock {\em IEEE Trans. Inform. Theory}, 33(6):788--801.

\bibitem[Davisson et~al., 1981]{Davisson-Longo-Sgarro-81}
Davisson, L.~D., Longo, G., and Sgarro, A. (1981).
\newblock The error exponent for the noiseless encoding of finite ergodic
  {M}arkov sources.
\newblock {\em IEEE Trans. Inform. Theory}, 27(4):431--438.

\bibitem[Dembo and Zeitouni, 1998]{Dembo-Zeitouni-98}
Dembo, A. and Zeitouni, O. (1998).
\newblock {\em Large deviations techniques and applications}, volume~38 of {\em
  Applications of Mathematics (New York)}.
\newblock Springer-Verlag, New York, second edition.

\bibitem[Dinwoodie, 1995]{Dinwoodie95}
Dinwoodie, I.~H. (1995).
\newblock A probability inequality for the occupation measure of a reversible
  {M}arkov chain.
\newblock {\em Ann. Appl. Probab.}, 5(1):37--43.

\bibitem[Donsker and Varadhan, 1975]{Donsker-Varadhan-I-75}
Donsker, M.~D. and Varadhan, S. R.~S. (1975).
\newblock Asymptotic evaluation of certain {M}arkov process expectations for
  large time. {I}. {II}.
\newblock {\em Comm. Pure Appl. Math.}, 28:1--47; ibid. 28 (1975), 279--301.

\bibitem[Ellis, 1984]{Ellis84}
Ellis, R.~S. (1984).
\newblock Large deviations for a general class of random vectors.
\newblock {\em Ann. Probab.}, 12(1):1--12.

\bibitem[Esscher, 1932]{Esscher32}
Esscher, F. (1932).
\newblock {On the probability function in the collective theory of risk}.
\newblock {\em Scandinavian Actuarial Journal}, 1932(3):175--195.

\bibitem[Fan et~al., 2018]{FJS18-Hoeffding}
Fan, J., Jiang, B., and Sun, Q. (2018).
\newblock {Hoeffding's lemma for Markov Chains and its applications to
  statistical learning}.

\bibitem[G\"{a}rtner, 1977]{Gartner77}
G\"{a}rtner, J. (1977).
\newblock On large deviations from an invariant measure.
\newblock {\em Teor. Verojatnost. i Primenen.}, 22(1):27--42.

\bibitem[Gillman, 1993]{Gillman93}
Gillman, D. (1993).
\newblock A {C}hernoff bound for random walks on expander graphs.
\newblock In {\em 34th {A}nnual {S}ymposium on {F}oundations of {C}omputer
  {S}cience ({P}alo {A}lto, {CA}, 1993)}, pages 680--691. IEEE Comput. Soc.
  Press, Los Alamitos, CA.

\bibitem[Hayashi and Watanabe, 2016]{HW16}
Hayashi, M. and Watanabe, S. (2016).
\newblock Information geometry approach to parameter estimation in {M}arkov
  chains.
\newblock {\em Ann. Statist.}, 44(4):1495--1535.

\bibitem[Hoeffding, 1963]{Hoeffding63}
Hoeffding, W. (1963).
\newblock Probability inequalities for sums of bounded random variables.
\newblock {\em J. Amer. Statist. Assoc.}, 58:13--30.

\bibitem[Horn and Johnson, 2013]{HJ13}
Horn, R.~A. and Johnson, C.~R. (2013).
\newblock {\em Matrix analysis}.
\newblock Cambridge University Press, Cambridge, second edition.

\bibitem[Jerrum et~al., 2001]{JSV01}
Jerrum, M., Sinclair, A., and Vigoda, E. (2001).
\newblock A polynomial-time approximation algorithm for the permanent of a
  matrix with non-negative entries.
\newblock In {\em Proceedings of the {T}hirty-{T}hird {A}nnual {ACM}
  {S}ymposium on {T}heory of {C}omputing}, pages 712--721. ACM, New York.

\bibitem[Jiang et~al., 2018]{FJS18-Bernstein}
Jiang, B., Sun, Q., and Fan, J. (2018).
\newblock {Bernstein's inequality for general Markov chains}.

\bibitem[Kahale, 1997]{Kahale97}
Kahale, N. (1997).
\newblock Large deviation bounds for {M}arkov chains.
\newblock {\em Combin. Probab. Comput.}, 6(4):465--474.

\bibitem[Kato, 1966]{Kato66}
Kato, T. (1966).
\newblock {\em Perturbation theory for linear operators}.
\newblock Die Grundlehren der mathematischen Wissenschaften, Band 132.
  Springer-Verlag New York, Inc., New York.

\bibitem[Kontoyiannis and Meyn, 2003]{KM03}
Kontoyiannis, I. and Meyn, S.~P. (2003).
\newblock Spectral theory and limit theorems for geometrically ergodic {M}arkov
  processes.
\newblock {\em Ann. Appl. Probab.}, 13(1):304--362.

\bibitem[Le\'{o}n and Perron, 2004]{LP04}
Le\'{o}n, C.~A. and Perron, F. (2004).
\newblock Optimal {H}oeffding bounds for discrete reversible {M}arkov chains.
\newblock {\em Ann. Appl. Probab.}, 14(2):958--970.

\bibitem[Lezaud, 1998]{Lezaud98}
Lezaud, P. (1998).
\newblock Chernoff-type bound for finite {M}arkov chains.
\newblock {\em Ann. Appl. Probab.}, 8(3):849--867.

\bibitem[Metropolis et~al., 1953]{MCMC53}
Metropolis, N., Rosenbluth, A.~W., Rosenbluth, M.~N., Teller, A.~H., and
  Teller, E. (1953).
\newblock Equation of state calculations by fast computing machines.
\newblock {\em The Journal of Chemical Physics}, 21(6):1087--1092.

\bibitem[Miasojedow, 2014]{Miaso14}
Miasojedow, B.~a. (2014).
\newblock Hoeffding's inequalities for geometrically ergodic {M}arkov chains on
  general state space.
\newblock {\em Statist. Probab. Lett.}, 87:115--120.

\bibitem[Miller, 1961]{Miller61}
Miller, H.~D. (1961).
\newblock A convexity property in the theory of random variables defined on a
  finite {M}arkov chain.
\newblock {\em Ann. Math. Statist.}, 32:1260--1270.

\bibitem[Moulos, 2019]{Moulos19-bandits-identification}
Moulos, V. (2019).
\newblock {Optimal Best Markovian Arm Identification with Fixed Confidence}.
\newblock In {\em 33rd Annual Conference on Neural Information Processing
  Systems}.

\bibitem[Nagaoka, 2005]{Nagaoka-05}
Nagaoka, H. (2005).
\newblock {The exponential family of Markov chains and its information
  geometry}.
\newblock In {\em Proceedings of The 28th Symposium on Information Theory and
  Its Applications (SITA2005)}, pages 1091--1095, Okinawa, Japan.

\bibitem[Paulin, 2015]{Paulin15}
Paulin, D. (2015).
\newblock Concentration inequalities for {M}arkov chains by {M}arton couplings
  and spectral methods.
\newblock {\em Electron. J. Probab.}, 20:no. 79, 32.

\bibitem[Rao, 2019]{Rao19}
Rao, S. (2019).
\newblock A {H}oeffding inequality for {M}arkov chains.
\newblock {\em Electron. Commun. Probab.}, 24:Paper No. 14, 11.

\bibitem[Watanabe and Hayashi, 2017]{WH17}
Watanabe, S. and Hayashi, M. (2017).
\newblock Finite-length analysis on tail probability for {M}arkov chain and
  application to simple hypothesis testing.
\newblock {\em Ann. Appl. Probab.}, 27(2):811--845.

\end{thebibliography}

\begin{appendices}
\section{Analyticity of Perron-Frobenius Eigenvalues and Eigenvectors}\label{app:anal}

Here we use the implicit function theorem in order to deduce
in~\autoref{lem:anal} that the Perron-Frobenius eigenvalue and eigenvectors are analytic functions of the entries of the matrix, at a level of generality adequate for our purposes.
\begin{lemma}\label{lem:anal}
Let $M \in \Realp{s \times s}$ be a nonnegative matrix
possessing the structure of~\autoref{lem:PF}, so in particular $M$ can be a nonnegative irreducible matrix.
\red 
Let $W$ range over $\mathbb{R}^{s \times s}$ in an open neighborhood of $M$.
\black 
Then $u_W, ~ \rho (W)$ and $v_W$ are analytic as functions of the entries of the matrix $W$ in 
\red an open neighborhood of \black $M$
where $W$ satisfies the conditions of~\autoref{lem:PF}.
\begin{proof}
We define the vector-valued function $F : \Real{(s+1)^2} \to \Real{2 (s + 1)}$
    \[
    F (W, u, \rho, v) :=
    \begin{bmatrix}
    (W\tran - \rho I) u \\
    \one\tran u - 1 \\
    (W - \rho I) v \\
    u\tran v - 1
    \end{bmatrix},
    \]
    where we use column vectors, and $\one$ denotes the all ones vector.
    At this point no assumptions are made about the structure of $W$.
    Note that each coordinate of the vector $F (W, u, \rho, v)$ is a multivariate polynomial of degree at most two, and hence each coordinate is an analytic function of $W, u, \rho$ and $v$.
    
    In addition $F (M, u_M, \rho (M), v_M) = 0$, and the Jacobian of $F$ with respect to $u, \rho, v$ evaluated at $W = M, u = u_M, \rho = \rho (M), v = v_M$ is
    \[
    J_{F, u, \rho , v} (M, u_M, \rho (M), v_M) =
    \begin{bmatrix}
    M\tran - \rho (M) I & - u_M & 0 \\
    \one\tran & 0 & 0 \\
    0 & - v_M & M - \rho (M) I \\
    v_M\tran & 0 & u_M\tran 
    \end{bmatrix}.
    \]
    We can easily verify that this Jacobian is left invertible.
    If $\begin{bmatrix} u\tran & \rho & v\tran \end{bmatrix}\tran$
    is in the kernel of $J_{F, u, \rho , v} (M, u_M, \rho (M), v_M)$,
    then $M\tran u = \rho (M) u + \rho u_M$, so if we multiply from the left with $v_M\tran$, we get that $\rho = 0$. In the same fashion, using~\autoref{lem:PF}, we can deduce that $u = v = 0$,
    and thus the kernel of the Jacobian is trivial.
    
    Then the analytic implicit function theorem guarantees that there exists
    a unique vector-valued function $g : \Real{s^2} \to \Real{2 s + 1}$ with each coordinate analytic, such that
    \[
    g (M) =
    \begin{bmatrix}
    u_M \\
    \rho (M) \\
    v_M
    \end{bmatrix}, ~\text{and}~
    F (W, g (W)) = 0, ~\text{for all $W$ in a neighborhood of}~ M.
    \]
    Let $k (M)$ denote the $1 \le k \le s$ corresponding to $M$ in the context of the discussion of the structure of $M$ in the statement of~\autoref{lem:PF}.
    Then, for $W$ in a sufficiently small neighborhood of $M$, the first $k (M)$ and the last $s+1$ coordinates of
     $g (W)$ have to be strictly positive.
     If we now restrict to those $W$ in such a neighborhood of $M$ that satisfy the conditions of~\autoref{lem:PF},
     then $k (W) \ge k (M)$ and, by~\autoref{lem:PF},
     $g (W)$ has to equal $\begin{bmatrix} u_W\tran & \rho (W) & v_W\tran \end{bmatrix}\tran$ for such matrices in this neighborhood of $M$.
\end{proof}
\end{lemma}
\section{Proofs from Section 2}
\label{app:misc-proofs}

\begin{proof}[Proof of~\autoref{lem:dual-map}]\hfill
\begin{enumerate}[label=(\alph*)]
    \item Each entry of $\tilde{P}_\theta$ is an analytic function of $\theta$, and the conclusion follows from~\autoref{lem:anal} in~\autoref{app:anal}.
    
    \item For any $x, y \in S$ such that $P (x, y) > 0$ we have
    \[
    \log P_\theta (x, y) = 
    \log P (x, y) + \theta f (y) - \Lambda (\theta) + \log v_\theta (y) - \log v_\theta (x).
    \]
    Differentiating with respect to $\theta$, and taking expectations with respect to $\pi_\theta \odot P_\theta$ we obtain
    \[
    \E_{(X, Y) \sim \pi_\theta \odot P_\theta} \frac{d}{d \theta} \log P_\theta (X, Y)
    = \pi_\theta (f) - \Lambda' (\theta).
    \]
   The conclusion follows because
    \[
    \E_{(X, Y) \sim \pi_\theta \odot P_\theta} \frac{d}{d \theta} \log P_\theta (X, Y)
    = \sum_x \pi_\theta (x) \frac{d}{d \theta} \left(\sum_y P_\theta (x, y)\right)
    = 0.
    \]
    
    \item For any $x, y \in S$ such that $P (x, y) > 0$ we have
    \[
    \frac{d^2}{d \theta^2} \log P_\theta (x, y) =
    - \Lambda'' (\theta) +
    \frac{d^2}{d \theta^2} \log v_\theta (y) -
    \frac{d^2}{d \theta^2} \log v_\theta (x).
    \]
    Taking expectations with respect to $\pi_\theta \odot P_\theta$ we obtain
    \begin{align*}
    \Lambda'' (\theta)
    &= - \E_{(X, Y) \sim \pi_\theta \odot P_\theta}\frac{d^2}{d \theta^2} \log P_\theta (X, Y) \\
    &= \E_{(X, Y) \sim \pi_\theta \odot P_\theta}\left(\frac{d}{d \theta} \log P_\theta (X, Y)\right)^2 \\
    &= \E_{(X, Y) \sim \pi_\theta \odot P_\theta}
    \left(
    f (Y) - \pi_\theta (f) +
     \frac{v_\theta (X)}{v_\theta (Y)}
    \frac{d}{d \theta}
    \frac{v_\theta (Y)}{v_\theta (X)}
    \right)^2.
    \end{align*}
    
    \item 
    Part (c) already ensures that $\Lambda (\theta)$ is convex. Moreover we see that
    \[
    \Lambda'' (\theta) = 0
    ~\text{for all}~ \theta \in (\theta_1, \theta_2),
    \quad\text{iff}\quad
    P_\theta = P_{\frac{\theta_1 + \theta_2}{2}}
    ~\text{for all}~ \theta \in (\theta_1, \theta_2).
    \]
    If such an interval $(\theta_1, \theta_2)$ exists, then
    we claim that we can enlarge it to the whole real line. To see this fix any $0 < \epsilon < \frac{\theta_2 - \theta_1}{2}$. Then using~\autoref{lem:comp}
    twice we obtain that for any $\theta \in (\theta_1, \theta_2)$
    \[
    P_{\theta \pm \epsilon} = 
    \left(P_\theta\right)_{\pm \epsilon} =
    \left(P_{\frac{\theta_1 + \theta_2}{2}}\right)_{\pm \epsilon} =
    P_{\frac{\theta_1 + \theta_2}{2} \pm \epsilon} =
    P_{\frac{\theta_1 + \theta_2}{2}}.
    \]
    By repeating this process we see that $P_\theta = P_0 = P$ for all $\theta \in \Real{}$.
    
    \red 
    Alternatively, if \black no such interval exists, then 
    $\Lambda' (\theta)$ is strictly increasing and $\Lambda (\theta)$ is strictly convex.
    \red Moreover, for \black $\theta_1 < \theta_2$ we have that
    $\pi_{\theta_1} (f) = \Lambda' (\theta_1) < \Lambda' (\theta_2) = \pi_{\theta_2} (f)$, and so $P_{\theta_1} \neq P_{\theta_2}$, establishing that in this case
    $\theta \mapsto P_\theta$ is an injection.

\end{enumerate}
\end{proof}

\begin{proof}[Proof of~\autoref{lem:KL-rep}]
\begin{align*}
\KL{\theta_1}{\theta_2}
&= \E_{(X,Y) \sim \pi_{\theta_1} \odot P_{\theta_1}}
\log \frac{P_{\theta_1} (X, Y)}{P_{\theta_2} (X, Y)} \\
&= \Lambda (\theta_2) - \Lambda (\theta_1) - \Lambda' (\theta_1) (\theta_2 - \theta_1) \\
&\quad + \E_{(X,Y) \sim \pi_{\theta_1} \odot P_{\theta_1}} 
\log \frac{v_{\theta_1} (Y)}{v_{\theta_1} (X)}
- \E_{(X,Y) \sim \pi_{\theta_1} \odot P_{\theta_1}} \log \frac{v_{\theta_2} (Y)}{v_{\theta_2} (X)} \\
&= \Lambda (\theta_2) - \Lambda (\theta_1) - \mu_1 (\theta_2 - \theta_1),
\end{align*}
where the second equality is using the calculations 
from the proof of~\autoref{lem:dual-map} (b).
\end{proof}

\begin{proof}[Proof of~\autoref{lem:conv-conj}]
From~\autoref{lem:dual-map} we have that
$\theta \mapsto \theta \mu - \Lambda (\theta)$ is either the linear function $\theta \mapsto (\mu - \pi (f)) \theta$, in which case the conclusion follows right away, or otherwise it is strictly concave.

In the latter case $\calM = (\mu_-, \mu_+)$ for some $\mu_- < \mu_+$.
If $\mu \in \calM$, then
$\theta = {\Lambda'}^{-1} (\mu)$ is the unique maximizer and the conclusion follows from~\autoref{lem:KL-rep}.
If $\mu = \mu_+$,
then the function keeps on growing as $\theta \to \infty$,
or equivalently as $\hat{\mu} \to \mu$, which in conjunction 
with the representation of the relative entropy rate from~\autoref{lem:KL-rep} establishes this case.
If $\mu > \mu_+$, then 
$\lim_{\theta \to \infty} \left(\theta \mu - \Lambda (\theta)\right) = 
\lim_{\theta \to \infty} \theta (\mu - \mu_+) + \lim_{\hat{\mu} \to \mu_+} 
\KL{\hat{\mu}}{\pi (f)} = \infty$.
The arguments are the same for the other two cases.
\end{proof}
\end{appendices}

\end{document}